\input amstex
\documentstyle{amsppt}
\magnification\magstep1
\NoRunningHeads
\hsize 6.5 truein
\vsize 8.0 truein
\define\k{\ref\key}
\define\pp{\pages}
\define\en{\endref}

\define\f{\frac}

\define\qq{\qquad}

\define\r1{\Bbb R^1}

\define\rn{\Bbb R^n}


\define\cm{\Cal M}


\define\nf{\infty}

\define\lp{ {L^p} }

\define\ra{\rightarrow}

\define\ga{\gamma}

\define\ve{\varepsilon}

\define\la{\lambda}

\topmatter
\title Best constants for uncentered maximal functions\endtitle
\author Loukas Grafakos$^*$ and Stephen Montgomery-Smith$^*$\endauthor
\affil University of Missouri, Columbia\endaffil
\address
Department of Mathematics, University of Missouri, Columbia, MO 65211
\endaddress
\email loukas\@msindy3.cs.missouri.edu, stephen\@mont.cs.missouri.edu
\endemail
\thanks $^*$Research partially supported by the NSF. \endthanks
\abstract
We  precisely evaluate  the operator norm of the
uncentered Hardy-Littlewood maximal function on $L^p(\Bbb R^1)$.
Consequently, we compute the operator norm of the
``strong'' maximal function  on $L^p(\Bbb R^n)$, and we observe that
the operator norm of the uncentered  Hardy-Littlewood maximal
function over balls on $L^p(\Bbb R^n)$ grows exponentially
 as $n \rightarrow \infty$.
\endabstract
\endtopmatter

\document

For a locally integrable function $f$
on $\rn$, let
$$
(\cm_n f)(x)= \sup_{B \ni x} \f{1}{|B|} \int_B |f(y)|\, dy,
$$
where the supremum is taken over all  closed balls $B$ that
contain the point $x$.
$\cm_n f $ is called the uncentered Hardy-Littlewood maximal
function of $f$ on $\rn$.
 In this paper we find the
precise value of the operator norm of $\cm_1$
on $L^p(\r1)$. It turns out that this operator norm is the
solution of  an  equation.
Our main result is the following:

\proclaim {Theorem}
For $1<p<\infty$, the operator norm of
$\cm_1$ : $\lp(\r1) \ra \lp(\r1)$ is the
unique positive solution of the   equation
$$
(p-1) \ x^p - p \ x^{p-1} -1 =0.\tag1
$$
\endproclaim

\bigskip

In order to prove our Theorem, we fix a nonnegative $f$ and we
introduce the left and
right maximal functions:
$$
(M_Lf)(x)= \sup_{a<x}\f{1}{x-a}\int_a^x f(t)\, dt \qq\text{and}\qq
(M_Rf)(x)= \sup_{b>x}\f{1}{b-x}\int_x^b f(t)\, dt.
$$
For the proof of the next result, known popularly as
the ``sunrise lemma'', we refer the reader to
Lemma (21.75) (i), Ch  VI in [HS].

\proclaim{Lemma 1}  Let $f\ge 0$ be in $L^1(\r1)$.  For each $\lambda > 0$, let
$C_\la  =\{x: (M_Lf)(x)>\la \}$\ and
$D_\la =\{x: (M_Rf)(x)>\la \}$.  Then
$$ \la |C_\la| = \int_{C_\la} f\, dt   \qq\text{and}\qq
   \la |D_\la| = \int_{D_\la} f\, dt. \tag2
$$
\endproclaim

\medskip

Now we are ready to prove the main lemma that leads to our Theorem.
This next result may be viewed as the ``correct'' weak type estimate
for the maximal function $\cm_1$.

\proclaim{Lemma 2}  Let  $f\ge 0$ be in  $L^1(\r1)$.  For each $\lambda > 0$,
let
$A_\la =\{x: (\cm_1f)(x)>\la \}$\ and
$B_\la =\{x: f(x)>\la \}$.
Then
$$
\la (|A_\la|+|B_\la|) \le \int_{A_\la}f \, dt+\int_{B_\la}f \, dt.\tag3
$$
\endproclaim

\noindent
To prove (3), first note that
$$
\sup (M_L,M_R)= \cm_1.\tag4
$$
For,
clearly $\sup (M_L,M_R)\le \cm_1$. On the other hand, it
is easy to see that for each real number $x$,\
$(\cm_1f)(x)$ is bounded by a convex combination of $(M_Lf)(x)$
and $(M_Rf)(x)$.

Now we add the two equalities in (2). Then using the fact that
$A_\la = C_\la \cup D_\la$ which follows from (4),  we obtain
$$
\la (|A_\la|+|C_\la \cap D_\la|) = \int_{A_\la}f\, dt +\int_{C_\la \cap
D_\la}f \, dt.\tag5
$$
Clearly $B_\la - (C_\la \cap D_\la)$
is a set of measure zero, and $f \le \la$ on
$(C_\la \cap D_\la)-B_\la$.  Therefore
$$
\int_{(C_\la \cap D_\la)-B_\la}f \, dt\le \la |(C_\la \cap D_\la)-B_\la|.
\tag6
$$
Equations (5) and (6) now imply equation $(3)$,
as required.

\bigskip

To prove the inequality in our Theorem, we require
the following fact.

\proclaim{Lemma 3} Let $f$\ and $g$\ be nonnegative functions on $\r1$.
Then if $p>1$, we have
$$ \int_0^\infty \lambda^{p-2} \int_{g(t) > \lambda} f(t) \, dt \, d\lambda
   =
   \frac1{p-1} \int_{\r1} f g^{p-1} \,dt ,$$
and if $p>0$, we have
$$ \int_0^\infty \lambda^{p-1} |\{g>\la\}| \, d\la
   =
   \frac1p \int_{\r1} g^p \,dt .$$
\endproclaim

\noindent
The first equality is easily proved, since by Fubini's theorem,
the left hand
side is
$$ \int_{-\infty}^\infty f(t) \int_0^{g(t)} \lambda^{p-2} \,d\lambda \, dt ,$$
which is readily seen to equal the right hand side.
The second equality is the special case of the first when $f=1$.

\bigskip

We now continue the proof of our Theorem.
Multiplying (3) by $\la^{(p-2)}$, integrating $\la$ from
$0$ to $\nf$, and applying Lemma 3, we obtain
$$
\f{1}{p}\|\cm_1f\|^p_p +\f{1}{p}\|f\|^p_p \le
\f{1}{p-1}\|f\|_p^p + \f{1}{p-1}\int_{\Bbb R^1}
 f(x) [(\cm_1f)(x)]^{p-1}\, dx,
$$
that is,
$$
(p-1)\|\cm_1f\|_p^p - \int_{\r1} f(x)[(\cm_1f)(x)]^{p-1} \, dx - \|f\|_p^p
\le 0.
$$
Applying H\" older's inequality with exponents $p$ and $p/(p-1)$ to
the second term, we obtain
$$
(p-1)\bigg(\f{\|\cm_1f\|_p}{\|f\|_p}\bigg)^{p}-
p\bigg(\f{\|\cm_1f\|_p}{\|f\|_p}\bigg)^{p-1}-1 \le 0, \tag7
$$
from which we conclude that $\f{\|\cm_1f\|_p}{\|f\|_p}\le c_p$,
where $c_p$ is
the unique positive solution of the (1).

To show that $c_p$ is in fact the operator norm of $\cm_1$ on
$\lp (\r1)$, we construct an example. Note that equality in (3)
is satisfied when $f$ is even symmetrically decreasing and equality
in (7) is satisfied when $\cm_1f$ is a multiple of $f$. We
are therefore led to the following example.
Let $f_{\ve, N}(t) = |t|^{-\f{1}{p}} \chi_{ \ve ,
 N}(|t|)$, where
$\chi_{\ve , N}$ is the characteristic
function of the  interval $[\ve ,N]$.
It can be easily seen that
$$
\lim_{\ve \ra 0}\lim_{N \ra \nf }\f{\|\cm_1f\|_p}{\|f\|_p} =
\cm_1(f_0)(1),\tag8
$$
where $f_0(t)=|t|^{-\f{1}{p}} \in L^1_{\text{loc}}$. An easy
calculation gives that
$$
\cm_1(f_0)(1)= \f{p}{p-1} \f{\ga^{\f{1}{p^\prime}}+1}{\ga +1},\tag9
$$
where $\ga $ is the unique positive solution of  the equation
$$
\f{p}{p-1} \f{\ga^{\f{1}{p^\prime}}+1}{\ga +1}= \ga^{-\f{1}{p}}.\tag10
$$
Using (9) and (10), it is a matter of simple
arithmetic to now show that $\cm_1(f_0)(1)$
is the unique positive root of  equation (1).
This completes the proof of our Theorem.

\bigskip

Before we conclude, we would like to make some remarks.
 Denote by  $x=(x_1, \dots , x_n)$ points in $\Bbb R^n$.
For a locally integrable function $f$ on $\rn$, define
$$
(\Cal N_n f)(x) = \sup\Sb a_1< x_1\\  b_1 > x_1\endSb \cdots
\sup\Sb a_n< x_n\\  b_n > x_n\endSb
{1\over (b_1-a_1)\cdots(b_n-a_n)} \int_{a_1}^{b_1} \cdots \int_{a_n}^{b_n}
f(y_1,\dots,y_n) \,dy_n\cdots dy_1 .
$$
$\Cal N_n$ is called the ``strong'' maximal
function on $\rn$. Clearly $\Cal N_1=\cm_1$. Observe that
$$ \Cal N_n \le \cm_1^{(1)} \circ \cdots \circ \cm_1^{(n)} ,$$
where $\cm_1^{(j)}$\ denotes the maximal operator $\cm_1$\ applied
to the $x_j$\ coordinate. This shows that the operator norm of
$\Cal N_n$\ on $L^p(\rn)$ is less than or equal to $c_p^n$.
By  considering the function
$$ g(x) = \prod_{j=1}^n f_{\epsilon,N}(x_j) ,$$
where  $f_{\epsilon,N}$ is as above, we obtain the converse
inequality.
We have therefore proved the following:

\proclaim{Corollary}
For $1<p<\infty$, the operator norm of
$\Cal N_n$ : $\lp(\rn) \ra \lp(\rn)$ is $c_p^n$, where $c_p$\ is the
unique positive solution of equation $(1)$.
\endproclaim

\bigskip

One can show that ${p \over p-1} < c_p < {2 p \over p-1}$. This implies
that the operator norm of $\Cal N_n$ on $L^p(\Bbb R^n)$ grows
exponentially with $n$, as $n \rightarrow \infty$.
Next, we observe that the same is true for the uncentered maximal
function $\cm_n$. There are several ways to see this. One way is by
considering the sequence of functions
$$
h_{\epsilon , N} (x)= |x|^{-\f{n}{p}} \chi_{ \epsilon  , N} (|x|) . $$
Let $U_n$ be the open unit ball in $\rn$.
For  $x \in \rn$, let $B_x=\frac{x}{2} +
\frac{|x|}{2}\overline{U_n}$. Then $x\in B_x$ and
$$
\big(\cm_n(h_{\ve , N}) \big)
(x)\geq   \frac{1}{|B_x|} \int_{B_x} |y|^{-\f{n}{p}}
\chi_{ \ve  , N}(|y|) \, dy =
\f{1}{|U_n|} \Bigl(\frac{2}{ |x|}\Bigr)^{n} \int_{B_x} |y|^{-\f {n}{p}}
\chi_{ \ve ,  N}(|y|) \, dy. \tag11
$$
Therefore for  $1<p<\infty$ and for all $\ve ,N>0$  we  have
$$ \align
\f{ \|\cm_n (h_{\ve , N})\|_{L^p}}{ \|
h_{\ve ,N}\|_{L^p}}
\ge  &\frac{2^n}{\|h_{\ve , N}\|_{L^p}|U_n| }
\Bigg\{ \int\limits_{r=0}^{\ +\nf } \int\limits_{S^{n-1}} \Big[
\frac{1}{r^n } \int\limits_{B_{r \phi}} |y|^{-\f{n}{p}} \chi_{
\epsilon ,  N}(|y|) dy \Big]^p d\phi  \  r^n \frac{dr}{r}
\Bigg\}^{\frac{1}{p}} \\
 =  &\frac{2^n}{\|h_{\epsilon , N}\|_{L^p} |U_n| }
\Bigg\{ \int\limits_{r=0}^{\ +\nf } \int\limits_{S^{n-1}} \Big[
\frac{1}{r^n } \int\limits_{t=0}^{ \ r} \int\limits_{S_\phi
(\frac{t}{r})}
 t^{-\f{n}{p}} \chi_{
\epsilon ,  N}(t) t^n \frac{dt}{t} d\theta \Big]^p d\phi \ r^n \frac{dr}{r}
\Bigg\}^{\frac{1}{p}}  , \tag12
\endalign$$
where $S_\phi(t)=\{ \theta \in S^{n-1} : |t\theta - \frac{\phi}{2} | \le
\frac{1}{2} \}$. By a change of variables (12)  is equal  to
$$\align
  &\frac{2^n}{\|h_{\ve , N}\|_{L^p} |U_n| }
\Bigg\{  \int\limits_{S^{n-1}}
 \int\limits_{r=0}^{\ +\nf }  \bigg[
 \int\limits_{t=0}^{ \ 1} \int\limits_{S_\phi
(t)}
  \chi_{
\ve ,  N}(rt) \ t^{\f{n}{p^\prime}}
 \frac{dt}{t} d\theta \bigg]^p   \frac{dr}{r}  d\phi
\Bigg\}^{\frac{1}{p}}\\
= &\f{2^n}{|U_n|}
\Bigg\{  \int\limits_{S^{n-1}}
\Bigg[\frac{\int\limits_{r=0}^{\ \nf }\big|(K_\phi * \chi_{\ve ,N})(r)
\big|^p
   \frac{dr}{r}}{ \int\limits_{r=0}^{\ \nf} \chi_{\ve ,N}^p(r) \f{dr}{r}
  }\Bigg] \frac{d\phi}{\omega_{n-1}}
\Bigg\}^{\frac{1}{p}}, \tag13
\endalign$$
where $K_\phi (t) = t^{n/p^\prime} \chi_{[0,1]}(t) \int_{S_\phi(t)}
|\theta|_B^{-n/p} d\theta$, $\omega_{n-1}=|S^{n-1}|=
\f{(n-1)\pi^{\f{n-1}{2}}}{\Gamma(\f{n+1}{2} )}$,
and $*$ denotes convolution on the multiplicative
group $ G=(\Bbb R^+, \frac{dt}{t}) $. If $K\ge 0$ on $G$,
the sequence of functions $\chi_{\epsilon ,
N}$ gives
equality in the convolution inequality $\|g*K\|_{L^p(G)} \le \|K\|_{L^1(G)}
\|g\|_{L^p(G)}$ as $\ve \ra 0 $ and $N \ra \nf$.
Therefore,  the
 expression inside brackets in (13) converges to $\|K_\phi\|_{L^1(G)}^p$
 as $\epsilon \ra 0 $ and $N \ra \nf$,
 and we obtain the estimate
$$\align
\lim \Sb \epsilon \rightarrow 0 \\ N \rightarrow \infty \endSb
\f{ \|\cm_n (h_{\ve , N})\|_{L^p}}{ \|
h_{\ve ,N}\|_{L^p}}
 &\ge \f{2^n}{|U_n|} \Bigg\{\int\limits_{S^{n-1}}
\bigg[ \int\limits_{0}^{\ 1} t^{\f{n}{p^\prime}}
\int\limits_{S_\phi(t)} \, d\theta \,  \f{dt}{t} \bigg]^p
\frac{d\phi}{\omega_{n-1}} \Bigg\}^{\frac{1}{p}}=
 \tfrac{n2^n}{\omega_{n-1}}\int\limits_{0}^{\ 1}
t^{\f {n}{p^\prime}} \int\limits \Sb S^{n-1} \\
\theta_1 \ge t\endSb d\theta  \, \f{dt}{t}              \\
 &= 2^{n}p^\prime
\tfrac{\omega_{n-2}}{
\omega_{n-1}}
 \int\limits_{0}^{\ 1} s^{\f{n}{p^\prime}} (1-s^2)^{\frac{n-3}{2}} \, ds =
2^{n-1} p^\prime \tfrac{\omega_{n-2}}{\omega_{n-1}}
  B(\tfrac{n}{2p^\prime}-\tfrac12 , \tfrac{n-3}{2}).\tag14
 \endalign$$
Stirling's formula
gives that expression
 (14)  is asymptotic to
$
  \Bigg\{{ \frac{ 4 ({\frac{1}{p^\prime}})^
{\frac{1}{p^\prime}} }{
(\frac{1}{p^\prime} +1)^
{(\frac{1}{p^\prime} +1)} }} \Bigg\}^{\frac{n}{2}}$
as $n \rightarrow \infty$
and since the number inside the braces above
is bigger than $1$ when
$1<p< \nf$, we also deduce exponential growth for
the operator norm of $\cm_n$ on $L^p(\Bbb R^n)$, as $n\rightarrow \infty$.

These remarks should be compared to the fact  that for $1<p<\infty$,
the operator norm of the
Hardy-Littlewood maximal function on $L^{p}(\rn)$ is
bounded above by some constant $A_{p}$ independent
of the dimension $n$\ (see [S] and [SS]).

\enddocument